\documentclass{article}
\usepackage{amsfonts,amssymb,latexsym,amsmath,amscd,euscript}
\newcounter{num}[section]
\newcommand{\Num}{\refstepcounter{num}%
\textbf{\arabic{section}.\arabic{num}}}

\newcommand{\Theorem}{\textbf{Theorem~}}
\newcommand{\Proof}{{\LARGE\emph{Proof}}}
\newcommand{\Def}{\textbf{Definition}}

\newcommand{\Lemma}{ \textbf{Lemma~}}
\newcommand{\Ex}{  \textbf{Example}}

\newcommand{\Prop}{\textbf{Proposition~}}
\newcommand{\Cor}{ \textbf{ Corollary~}}

\newcommand{\GL}{{\mathrm{GL}}}
\newcommand{\SL}{{\mathrm{SL}}}

\newcommand{\UT}{{\mathrm{UT}}}

\newcommand{\Xb}{{\Bbb X}}
\newcommand{\Yb}{{\Bbb Y}}
\newcommand{\Sb}{{\Bbb S}}

\newcommand{\Ib}{{\Bbb I}}
\newcommand{\Jb}{{\Bbb J}}
\newcommand{\Pb}{{\Bbb P}}
\newcommand{\Gb}{{\Bbb G}}

\newcommand{\al}{\alpha}

\newcommand{\Om}{\mathrm{O}}
\newcommand{\Spm}{\mathrm{Sp}}

\renewcommand{\leq}{\leqslant}
\renewcommand{\geq}{\geqslant}

\newcommand{\Jc}{{\mathcal J}}
\newcommand{\Ac}{{\mathcal A}}
\newcommand{\Fc}{{\mathcal F}}
\newcommand{\Nc}{{\mathcal N}}
\newcommand{\Xc}{{\mathcal X}}
\newcommand{\Mc}{{\mathcal M}}
\newcommand{\Pc}{{\mathcal P}}
\newcommand{\Ax}{{\mathfrak A}}

\newcommand{\Mat}{{\mathrm{Mat}}}
\newcommand{\diag}{{\mathrm{diag}}}

\begin{document}
\Large

\title{Fields of invariants for unipotent radicals of  parabolic subgroups}
\author{A. N. Panov\footnote{The work is supported by the RFBR grant   20-01-00091a}}
\date{}
 \maketitle

\begin{abstract}
	
	The paper is devoted to the problem of finding free  generators in the fields of invariants  for actions of unipotent groups on affine  varieties. We consider the case when the unipotent group is  the unipotent  radical in an arbitrary parabolic subgroup in the reductive group of  classical type 
	$\GL(n)$,~ $\SL(n)$, ~ $\Om(n)$ or $\Spm(2n)$. In the explicit form, we present a system of free generators in the field of invariants for the action of the unipotent radical on the reductive group by conjugation.  
\end{abstract}

\section{Introduction}

Let  $G$ be a split reductive group defined over a field $K$ of characteristic zero. Let  $U$ be  the unipotent radical of an arbitrary parabolic subgroup  $P$ in $G$.  The group  $U$ acts on $G$ by conjugation. There is a representation of the group  $U$ on the space  $\Ac=K[G]$   by the formula $\rho(g)f(x) = f(g^{-1}xg)$,~ $g\in U$, ~ $x\in G$.
This representation extends to the action of  $U$ on the field  $\Fc=K(G)$ of rational functions on $G$.

In the present paper we  consider the cases when  $G$  is one of the following groups $\GL(n)$, ~$\SL(n)$, ~$\Om(n)$ or $\Spm(2n)$. We aim to construct a system of free generators of the field of invariants  $\Fc^U$ in an explicit form as polynomials in matrix entries.

In the second section of this  paper, we solve this problem  for the case of the general linear  group  $\GL(n)$. In  Theorem  \ref{theoremone}, we construct the system of free generators  $\{\Jc_{ij}:~ (i,j)\in\Sb\}$ for the field $\Fc^U$, where   $\Jc_{ij}$ are determinants of special form. We obtain the similar system of generators for  $\SL(n)$ (see Corollary  \ref{corone}). 
The case when $U$ is the unitriangular group (i.e. the group of of upper triangular matrices with ones on the diagonal) is treated by the author in the paper \cite{PV1}.

In the third section, we consider the case of  $G=\Om(n)$ or $G=\Spm(2n)$. In Theorem    \ref{theoremtwo}, using the restrictions $\{\Jc_{ij}^\circ\}$ of polynomials $\{\Jc_{ij}\}$  on $G$, we construct the system of free generators for the  field of $U$-invariants on the group $G$.

Since   the group $U$ is  unipotent,  the field of invariants $\Fc^U$ is   a field of fractions of the algebra  of invariant regular functions $\Ac^U$ ~\cite[Theorem 3.3]{PV}. The algebra  $\Ac^U$ is finitely generated~ \cite[Theorem 3.13]{PV}.  The problem of description of generators for the  algebra $\Ac^U$  is still an unsolved problem. 

In conclusion of Introduction we emphasize a few papers which are devoted to the theory of invariants of unipotent groups \cite{PV2, Mi,Pom,Panyu}.

\section{Fields of invariants for  radicals of  parabolic subgroups in  $\GL(n)$}

Let  $K$ be a field of  characteristic zero and   $P$ be a parabolic subgroup of the general linear group  $G=\GL(n)$. There is a decomposition   $P=LU$ of the group $P$ into a semidirect product of its unipotent radical  $U$ and a Levi subgroup  $L$. 
Our aim is to construct a system of free generators of the field $U$-invariants $\Fc^{U}$ explicitly.

The parabolic subgroup  $P$ is defined by a partition of  the integer segment   $[1,n]$   into a system of  consecutive segments $$[1,n] =I_1\sqcup I_2\sqcup\ldots\sqcup I_\ell.$$    Put  $n_i=|I_i|$ for each $1\leq i\leq \ell$.
The linear space $\Mc=\Mat(n)$ of all $(n\times n)$ matrices  is a direct sum  $$\Mc=\bigoplus_{1\leq k,m\leq \ell} \Mc_{km}, $$ where $\Mc_{km}$ is spanned by  the matrix units $E_{ij}$,~ $(i,j)\in I_k\times I_m$. 
The Levi subgroup $L$ is isomorphic to the direct  product 
$$\GL(n_1)\times\ldots \times \GL(n_\ell).$$ 
The unipotent subgroup $U$ consists of matrices $E+A$, where $E$  is the identity matrix and  $A$ belongs to the sum of subspaces  $\Mc_{km}$, ~  $k<m$.

For each $i\in [1,n]$, let  $i'=n+1-i$ be the symmetric number to $i$ with respect to the center of the segment. Respectively, for any  $T\subseteq [1,n]$, we define $T'=\{i':~ i\in  T\}$.

 We consider the set $\Sb$ of pairs $(i,j)$  such that  $$i\in I_k,~~ j\in I_m',~~ k\geq m.$$ Let  $S$ be the intersection of the group $\GL(n)$ with the subspace in $\Mc$  spanned by  $E_{ij}$, ~ $(i,j)\in \Sb$.\\
 \\
 \Prop\Num\label{open}. The subset   \begin{equation}\label{gSg}
 \bigcup_{g\in U} gSg^{-1}
 \end{equation}
 is dense in  $\GL(n)$.\\
 \Proof.  Let $\Nc=\UT(n)$,~ $\Nc_L=\Nc\cap L$ and  $B$ be the subgroup of upper triangular matrices in $\GL(n)$. The Bruhat cell  $\Nc w_0B$, where   $w_0$ is  the element of greatest length in the Weyl group,  is dense in  $\GL(n)$.  The subgroup  $\Nc$ decomposes into a product  $\Nc=U\cdot \Nc_L$. Then  
 $\Nc w_0B = U\left(\Nc_Lw_0 B\right)$.  By direct calculations one can verify  $\Nc_Lw_0 B\subset S$ (see Example \ref{exone}).  Then  $\Nc w_0B\subset US$.  
 Since $S$ is stable with respect to the right multiplication by  $U$, the subset  $US$ coincides with the subset (\ref{gSg}).  Therefore (\ref{gSg}) is dense in  $\GL(n)$.    ~$\Box$

 Denote by  $\pi$  the restriction map  $\Ac\to K[S]$.  The map $\pi$ establishes an  embedding  $\Ac^U\to K[S]$ and it is extended to an embedding  $\pi: \Fc^U\to K(S)$.

Let  $\{x_{ij}\}$ be the system of standard coordinate functions on  $\Mc$.
Define the matrix  $\Xb=(x_{ij})_{i,j=1}^n$. Consider the adjugate  matrix  $$\Xb^* = (x^*_{ij})_{i,j=1}^n.$$ 
We have $\Xb\Xb^*=\Xb^*\Xb=\det(\Xb)E.$ 

Define the order relation  on the set of pairs $\Sb$ as follows:
\begin{equation}\label{Order}
(i_1,j_1)\prec (i,j)~~\mbox{if}~~ j_1<j, ~~\mbox{or}~~ j_1=j~~\mbox{and}~~ i_1>i.
\end{equation}
Denote  $s_{ij} =\pi(x_{ij})$ for $(i,j)\in \Sb$.

Our aim is to attach the  $U$-invariant $\Jc_{ij}$ in $K[\Mc]$ to each  pair  $(i,j)\in \Sb$. Decompose $\Sb$ into two subsets $\Sb=\Sb_0\sqcup \Sb_1$, where
$\Sb_0$ consists of all $(i,j)\in\Sb$, which is lying   on or above the anti-diagonal, and
 $\Sb_1$ consists of all $(i,j)\in\Sb$, which is lying strictly below  the anti-diagonal.

Let  $(i,j)\in\Sb_0$.  Then  $(i,j)$  is lying on or above the anti-diagonal. In this case $i+j\leq n+1$ (equivalent to  $i\leq j'$). 
Since $i\leq j'$,  the number  $i$ is less than any number from the segment $[j'+1, n]$. Observe that  $[1,j-1]'=[j'+1, n]$. Consider the system of $j$ rows $\{i\}\sqcup [j'+1, n]$. 

For each   $(i,j)\in\Sb_0$, we define the polynomial   $\Jc_{ij}$  as the minor of the matrix  $\Xb$ defined by the systems of columns  $[1,j]$ and rows $\{i\}\cup [j'+1,n]$. 

 Let $(i,j)\in\Sb_1$. Then $(i,j)$ is lying below the anti-diagonal. In this case  $i+j > n+1$ (equivalent to  $j>i'$ ). Then $j=i'+ (j-i')$.
  Consider the  $(j\times j)$-matrix  $$\Yb_{ij}= \left(\begin{array}{c} \Xc_{i',j}\\ ----- \\
\Xc^*_{j-i', j}\end{array}\right),$$ where  $\Xc_{i',j}$ is  the submatrix   $\Xb$ associated with  the columns $[1,j]$  and the last  $i'$ rows, 
   $\Xc^*_{j-i',j}$ is  the submatrix of  $\Xb^*$ associated  with  the columns  $[1,j]$ and the last $j-i'$ rows.  
For each $(i,j)\in\Sb_1$, we define 
$$\Jc_{ij} = \det \Yb_{ij}.$$
\Prop\Num\label{Jinv}. The polynomials  $\{\Jc_{ij}:~~(i,j)\in \Sb\}$ are  $U$-invariants.
\\
\Proof. Analogously to    \cite[Proposition 1]{PV1}. 
\\
\Prop\Num\label{Jnonzero}.~ $\Jc_{ij}\ne 0$ for each $(i,j)\in \Sb$. \\
\Proof. The statement is obvious for   $(i,j)\in\Sb_0$. Let  $(i,j)\in\Sb_1$. Consider the matrix  $A=(a_{st})$, where  
$a_{st}=1$ if  $(s,t)$ lies on the anti-diagonal or   $(s,t) = (i+k,j-k)$ for $ k\in [0,n-i]$, otherwise  $a_{st}=0$. Then in the determinant $\Jc_{ij}(A)$
all below the anti-diagonal entries of
 are zero  and the anti-diagonal entries are non-zero. Therefore $\Jc_{ij}(A)\ne 0$.~ 
 For instance, in Example \ref{exone}, for $i=j=4$  we have 
 
{\small $$A=\left(\begin{array}{ccccc}
 	0&0&0&0&1\\
 	0&0&0&1&0\\
 	0&0&1&0&0\\
 	0&1&0&1&0\\
 		1&0&1&0&0\\
 \end{array}\right),\qquad A^*=\left(\begin{array}{ccccc}
 \times&\times&\times&\times&1\\
 \times&\times&\times&1&0\\
 \times&\times&1&0&0\\
 \times&1&0&0&0\\
 1&0&0&0&0\\
 \end{array}\right), \qquad
 \Jc_{44}(A)=\left|\begin{array}{cccc}
 	0&1&0&1\\
 1&0&1&0\\
 \times&1&0&0\\
 1&0&0&0\\
 \end{array}\right| = 1\ne 0.~~\Box$$}
\Cor\Num. ~ $\pi(\Jc_{ij})\ne 0$ for each $(i,j)\in \Sb$. \\
The proof follows from the fact that $\pi$ is an embedding of  $\Fc^U$ to $K(S)$.\\
\Lemma\Num\label{minorstar}. Let   $M_{I,J}(\Xb^*)$ be the minor of $\Xb^*$ with systems of rows $I$ and columns $J$. Then  if $I$ is a segment and it ends in  $n$, then  $M_{I,J}(\Xb^*)$  is invariant under  the right multiplication of  $\Xb$ by $\UT(n)$.\\
\Proof. If $I$ is a segment and it ends in  $n$, then for any  $(n\times n)$-matrix $\Yb$  and $g\in \UT(n)$, we have  $M_{I,J} (g\Yb) = M_{I,J} (\Yb)$. Therefore
$$M_{I,J} \left((\Xb g)^*\right) = M_{I,J} (g^*\Xb^*) = M_{I,J} (\Xb^*).~~\Box$$

The subset  $S$ is stable under the right multiplication by  $\UT(n)$.
Denote by  $K(S)^{\UT(n)}$ the subfield of $\UT(n)$-invariants in   $K(S)$. \\
\Cor\Num\label{corrstar}.   If $I$ is a segment and it ends in  $n$, then  $\pi(M_{I,J}(\Xb^*))$ belongs to  $K(S)^{\UT(n)}$. \\
\Def~\Num. Let  $\{\xi_\al: \al\in \Ax\}$  and   $\{\eta_\al: \al\in \Ax\}$ be two finite systems of free generators of an extension  $F$ of the field  $K$. Let $\prec$ be a linear order relation on  $\Ax$. We say that the second system of generators is obtained from the first one by a triangular transformation if  each  $\eta_\al$ can be presented in the form 
\begin{equation}\label{etaxiorder}
\eta_\al=\phi_\al\xi_\al+\psi_\al,
\end{equation}
where $\phi_\al\ne 0$ and $\phi_\al$,~ $\psi_\al$ belong to the subfield generated by  $\{\xi_\beta: ~ \beta\prec \al\}$.

Denote by $\Pc_0$ the subfield of $\Fc^U$ generated by   $\{\Jc_{ij}:~~ (i,j)\in \Sb_0\}$.  \\
\\ 
\Prop\Num\label{fzero}. \\
1) The subfield  $K(S)^{\UT(n)}$ is freely generated by the elements 
 $$\{\pi(\Jc_{ij}):~~ (i,j)\in \Sb_0\}.$$
 2) The map  $\pi$ establishes  an isomorphism of the subfield  $\Pc_0$ of $\Fc^U$ to  $K(S)^{\UT(n)}$.\\
\Proof. Denote by  $S_0$ the subset of  $S$ that  consists of  matrices  $(s_{ij})$ with  $s_{ij}=0$ for  all  $(i,j)\in \Sb_1$. 
The subset $S_0\cdot\UT(n)$ is dense in $S$. The restriction  map  $\pi_0$  from  $K[S]$ to  $K[S_0]$ defines  an embedding of the subfield  $K(S)^{\UT(n)}$ into the field $K(S_0)$.   Each polynomial in  $\{\pi(\Jc_{ij}):~ (i,j)\in \Sb_0\}$ is a $\UT(n)$-invariant and its restriction $ \pi_0\pi(\Jc_{ij})$ on  $S_0$  equals to 
$$\pm s_{n1}s_{n-1,2}\ldots s_{j'+1,j-1}s_{ij}.$$ 
For instance, in Example \ref{exone}  
{\small $$S= \left\{\left(\begin{array}{ccccc}
	0&0&0&0&s_{15}\\
	0&0&s_{23}&s_{24}&s_{25 }\\
	0&0&s_{33}&s_{34}&s_{35 }\\
	s_{41}&s_{42}&s_{43}&s_{44}&s_{45 }\\
	s_{51}&s_{52}&s_{53}&s_{54}&s_{55 }\\
\end{array}\right)\right\},\qquad S_0=\left\{\left(\begin{array}{ccccc}
0&0&0&0&s_{15}\\
0&0&s_{23}&s_{24}&0\\
0&0&s_{33}&0&0\\
s_{41}&s_{42}&0&0&0\\
s_{51}&0&0&0&0\\
\end{array}\right)\right\}.$$}
In the case  $i=2$,~$j=3$, we obtain 
 $$\pi_0\pi(\Jc_{23})=\left|\begin{array}{ccc}
	0&0&s_{23}\\
	s_{41}&s_{42}&0\\
	s_{51}&0&0\\
	\end{array}\right|= -s_{51} s_{42} s_{23}.$$
The system of elements  $\{\pi_0\pi(\Jc_{ij}):~ (i,j)\in \Sb_0\}$ is obtained by a triangular transformation with respect to the order relation  $\prec$ (see (\ref{Order})) from  the system  
$\{s_{ij}:~ (i,j)\in\Sb_0\}$ of generators of the field  $K(S_0)$. This implies  $\pi_0$ is an isomorphism of  $K(S)^{\UT(n)}$ to the field  $K(S_0)$. ~ $\Box$\\
\\
\Theorem\Num\label{theoremone}. The system of polynomials  $\{\Jc_{ij}: ~(i,j)\in \Sb\}$ freely generates over $K$ the field of invariants  $\Fc^{U}$ for the group $\GL(n)$.\\
\Proof. \\
\textit{Item 1}. Let us show that the field  $K(S)$ is freely  generated  over $K(S)^{\UT(n)}$ 
by the system of elements $\{s_{ij}:~ (i,j)\in\Sb_1\}$. 

For each  $(i,j)\in \Sb_0$, we expand the minor  $\Jc_{ij}$ along its first row.
We get  $\Jc_{ij} =\pm A_{ij} \cdot  x_{ij}+\Psi_{ij}$.  The cofactor  $A_{ij}$ and   $\Psi_{ij}$, which is the sum of the other cofactors multiplied by the respective entries, both are  expressions in the $x_{st}$, ~$(s,t)\prec (i,j)$. Therefore, the system of elements 
$$\{\pi(\Jc_{ij}):~ (i,j)\in \Sb_0\}\sqcup \{s_{ij}:~ (i,j)\in\Sb_1\}$$
is obtained by a triangular transformation from the system of free generators $\{s_{ij}:~ (i,j)\in\Sb\}$ of the field  $K(S)$. This proves the statement of Item 1.\\
\textit{Item 2}. Let us show that the field $\Fc^U$ is freely generated over the subfield  $\Pc_0$ by the system of elements    $\{\Jc_{ij}:~ (i,j)\in\Sb_1\}$.
It is sufficient to prove that the system of elements   $\{\pi(\Jc_{ij}):~ (i,j)\in\Sb_1\}$
 is obtained from  $\{s_{ij}:~(i,j)\in\Sb_1\}$ by a triangular transformation over the field  $K(S)^{\UT(n)}$. We prove using the induction method  that the formula of type  (\ref{etaxiorder}) is valid for the order relation  (\ref{Order}).
 
 For the smallest element $(n,1)$, we have  $\Jc_{ij}=x_{n1}$, and $\pi(\Jc_{ij})=s_{n1}$ belongs to  $ K(S)^{\UT(n)}$. 
 
 Assume that the statement is true for all elements  $\prec (i,j)$. Let us prove for  $(i,j)\in\Sb_1$. 
Expand the determinant  $\Jc_{ij}$ along its first row.
We get 
$$\Jc_{ij} =\pm M \cdot  x_{ij}+\Psi_{ij}, $$
where $M$ is the cofactor  for  $(i,j)$,  and   $\Psi_{ij}$ is the sum of the other entries of the first row multiplying to their cofactors. For instance, 
in Example \ref{exone} for $i=5$,~ $j=3$ we have 
{\small  
$$\Jc_{53}=\left|\begin{array}{ccc}
x_{51}&x_{52}&x_{53}\\
x^*_{41}&x^*_{42}&x^*_{43}\\
x^*_{51}&x^*_{52}&x^*_{53}\\
\end{array}\right| = M x_{53} + \Psi_{53},$$
$$M= \left|\begin{array}{cc}
x^*_{41}&x^*_{42}\\
x^*_{51}&x^*_{52}\\
\end{array}\right|,\qquad \Psi_{53} = - \left|\begin{array}{cc}
x^*_{41}&x^*_{43}\\
x^*_{51}&x^*_{53}\\
\end{array}\right| x_{52} + \left|\begin{array}{cc}
x^*_{42}&x^*_{43}\\
x^*_{52}&x^*_{53}\\
\end{array}\right| x_{51}.$$}

Restricting the above equality to  $S$,  we obtain
\begin{equation}\label{pijc}
\pi(\Jc_{ij}) =\pm \pi(M)\cdot  s_{ij}+\pi(\Psi_{ij}).
\end{equation} 

 In the case $i<n$, we have $M=\Jc_{i+1,j-1}$ and $M\ne 0$. 
According to induction assumption, $\pi(\Jc_{i+1,j-1})$  is an expression over the field  $K(S)^{\UT(n)}$  of   $s_{kt}$, ~ $(k,t)\in\Sb_1$, ~ $(k,t)\prec (i,j)$.
 
  In the case $i=n$, the polynomial $M$ coincides with the minor 
$ M_{I,J}(\Xb^*)$ for $J=[1,j-1]$ and $I=J'=[j'+1,n]$. The segment  $I$ ends in $n$, therefore $ \pi(M_{I,J}(\Xb^*))\in K(S)^{\UT(n)}$ (see Corollary  \ref{corrstar}). 

The sum  $\Psi_{ij}$ is an expression of  $x_{kt}$, ~$(k,t)\prec (i,j)$, and   minors of the form  $M_{I',J'}(\Xb^*)$, where the segments of rows  $I'$  end in  $n$. By the induction assumption and Corollary  \ref{corrstar}, we verify that  $\pi(\Psi_{ij})$  is an expression over the field  $K(S)^{\UT(n)}$ by  $s_{kt}$, ~ $(k,t)\in\Sb_1$, ~ $(k,t)\prec (i,j)$. The statement of Item 2 is proved.  ~$\Box$
 \\
 \Cor\Num\label{corone}. The system of polynomials  $\{\Jc_{ij}: ~(i,j)\in \Sb,~ (i,j)\ne (1,n)\}$ freely generates over $K$ the field of invariants  $\Fc^{U}$ for the group $\SL(n)$.
 \\
 \\
\Ex~\Num\label{exone}. Let  $n=5$ and the parabolic subgroup $P$ be defined by the partition $[1,5] = \{1\}\sqcup \{2,3\}\sqcup \{4,5\}$. In the  $(5\times 5)$-table $\Pb$, by  the symbol 
$"\times"\,$ (respectively, by the symbol $"*"\,$), we mark the cells that correspond to  the subgroup $L$ (respectively, to the subgroup $U$).  
{\small 
$$\Pb =\left(\begin{array}{ccccc}
\times&*&*&*&*\\
0&\times&\times&*&*\\
0&\times&\times&*&*\\
0&0&0&\times&\times\\
0&0&0&\times&\times\\
\end{array}\right),
~~ 
S= \left\{\left(\begin{array}{ccccc}
0&0&0&0&s_{15}\\
0&0&s_{23}&s_{24}&s_{25 }\\
0&0&s_{33}&s_{34}&s_{35 }\\
s_{41}&s_{42}&s_{43}&s_{44}&s_{45 }\\
s_{51}&s_{52}&s_{53}&s_{54}&s_{55 }\\
\end{array}\right)\right\},
~~ 
S^*= \left\{\left(\begin{array}{ccccc}
s^*_{11}&s^*_{12}&s^*_{13}&s^*_{14}&s^*_{15 }\\
s^*_{21}&s^*_{22}&s^*_{23}&s^*_{24}&s^*_{25 }\\
s^*_{31}&s^*_{32}&s^*_{33}&0&0\\
s^*_{41}&s^*_{42}&s^*_{43}&0&0\\
s^*_{51}&0&0&0&0\\
\end{array}\right)\right\}, $$}

{\small 
	$$\Nc_L= \left(\begin{array}{ccccc}
	1&0&0&0&0\\
	0&1&\times&0&0\\
	0&0&1&0&0\\
	0&0&0&1&\times\\
	0&0&0&0&1\\
	\end{array}\right),
	~~ 
	\Nc_L w_0 B = \left(\begin{array}{ccccc}
	1&0&0&0&0\\
	0&1&\times&0&0\\
	0&0&1&0&0\\
	0&0&0&1&\times\\
	0&0&0&0&1\\
	\end{array}\right) 	
\left(\begin{array}{ccccc}
0&0&0&0&\times\\
0&0&0&\times&\times\\
0&0&\times&\times&\times\\
0&\times&\times&\times&\times\\
\times&\times&\times&\times&\times\\
\end{array}\right) \subset  S. $$}

The generators  $\Jc_{ij}$,~ $(i,j)\in\Sb$, and their restrictions  $\pi(\Jc_{i,j})$  on $S$ have the following form:\\
{\small 
 $\Jc_{51}=x_{51}$,~~  $\pi(\Jc_{51})=s_{51}$,\qquad\qquad $\Jc_{41}=x_{41}$, ~ $\pi(\Jc_{41})=s_{41}$,\\
$\Jc_{52}=\left|\begin{array}{cc}
x_{51}&x_{52}\\
x^*_{51}&x^*_{52}\\
\end{array}\right|$, ~ $\pi(\Jc_{52})=\left|\begin{array}{cc}
s_{51}&s_{52}\\
s^*_{51}&0\\
\end{array}\right|$,\qquad\qquad
$\Jc_{42}=\left|\begin{array}{cc}
x_{41}&x_{42}\\
x_{51}&x_{52}\\
\end{array}\right|$,  ~~ $\pi(\Jc_{42})=\left|\begin{array}{cc}
	s_{41}&s_{42}\\
	s_{51}&s_{52}\\
\end{array}\right|$, \\
 $\Jc_{53}=\left|\begin{array}{ccc}
x_{51}&x_{52}&x_{53}\\
x^*_{41}&x^*_{42}&x^*_{43}\\
x^*_{51}&x^*_{52}&x^*_{53}\\
\end{array}\right|$, ~~ $\pi(\Jc_{53})=\left|\begin{array}{ccc}
s_{51}&s_{52}&s_{53}\\
s^*_{41}&s^*_{42}&s^*_{43}\\
s^*_{51}&0&0\\
\end{array}\right|$,
\\
$\Jc_{43}=\left|\begin{array}{ccc}
x_{41}&x_{42}&x_{43}\\
x_{51}&x_{52}&x_{53}\\
x^*_{51}&x^*_{62}&x^*_{53}\\
\end{array}\right|$, ~ $\pi(\Jc_{43})=\left|\begin{array}{ccc}
s_{41}&s_{42}&s_{43}\\
s_{51}&s_{52}&s_{53}\\
s^*_{51}&0&0\\
\end{array}\right|$;\qquad\qquad $\Jc_{33}=\left|\begin{array}{ccc}
x_{31}&x_{32}&x_{33}\\
x_{41}&x_{42}&x_{43}\\
x_{51}&x_{52}&x_{53}\\
\end{array}\right|$,~~ $\pi(\Jc_{33})=\left|\begin{array}{ccc}
0&0&s_{33}\\
s_{41}&s_{42}&s_{43}\\
s_{51}&s_{52}&s_{53}\\
\end{array}\right|$,\\
\\
$\Jc_{23}=\left|\begin{array}{ccc}
x_{21}&x_{22}&x_{23}\\
x_{41}&x_{42}&x_{43}\\
x_{51}&x_{52}&x_{53}\\
\end{array}\right|$,~~ $\pi(\Jc_{23})=\left|\begin{array}{ccc}
0&0&s_{23}\\
s_{41}&s_{42}&s_{43}\\
s_{51}&s_{52}&s_{53}\\
\end{array}\right|$,\\
\\
$\Jc_{54}=\left|\begin{array}{cccc}
x_{51}&x_{52}&x_{53}&x_{54}\\
x^*_{31}&x^*_{32}&x^*_{33}&x^*_{34}\\
x^*_{41}&x^*_{42}&x^*_{43}&x^*_{44}\\
x^*_{51}&x^*_{52}&x^*_{53}&x^*_{54}\\
\end{array}\right|$, ~  $\pi(\Jc_{54})=\left|\begin{array}{cccc}
s_{51}&s_{52}&s_{53}&s_{54}\\
s^*_{31}&s^*_{32}&s^*_{33}&0\\
s^*_{41}&s^*_{42}&s^*_{43}&0\\
s^*_{51}&0&0&0\\
\end{array}\right|$,\qquad\qquad \\
\\
$\Jc_{44}=\left|\begin{array}{cccc}
x_{41}&x_{42}&x_{43}&x_{44}\\
x_{51}&x_{52}&x_{53}&x_{54}\\
x^*_{41}&x^*_{42}&x^*_{43}&x^*_{44}\\
x^*_{51}&x^*_{52}&x^*_{53}&x^*_{54}\\
\end{array}\right|$, ~  $\pi(\Jc_{44})=\left|\begin{array}{cccc}
s_{41}&s_{42}&s_{43}&s_{44}\\
s_{51}&s_{52}&s_{53}&s_{54}\\
s^*_{41}&s^*_{42}&s^*_{43}&0\\
s^*_{51}&0&0&0\\
\end{array}\right|$,
\\
\\
$\Jc_{34}=\left|\begin{array}{cccc}
x_{31}&x_{32}&x_{33}&x_{34}\\
x_{41}&x_{42}&x_{43}&x_{44}\\
x_{51}&x_{52}&x_{53}&x_{54}\\
x^*_{51}&x^*_{52}&x^*_{53}&x^*_{54}\\
\end{array}\right|$, ~  $\pi(\Jc_{34})=\left|\begin{array}{cccc}
0&0&s_{33}&s_{34}\\
s_{41}&s_{42}&s_{43}&s_{44}\\
s_{51}&s_{52}&s_{53}&s_{54}\\
s^*_{51}&0&0&0\\
\end{array}\right|$,\\
\\  
$\Jc_{24}=\left|\begin{array}{cccc}
x_{21}&x_{22}&x_{23}&x_{24}\\
x_{31}&x_{32}&x_{33}&x_{34}\\
x_{41}&x_{42}&x_{43}&x_{44}\\
x_{51}&x_{52}&x_{53}&x_{54}\\
\end{array}\right|$, ~ 
 $\pi(\Jc_{24})=\left|\begin{array}{cccc}
  0&0&s_{23}&s_{24}\\
 0&0&s_{33}&s_{34}\\
s_{41}&s_{42}&s_{43}&s_{44}\\
s_{51}&s_{52}&s_{53}&s_{54}\\
\end{array}\right|$, \\
\\
$\Jc_{55} = \left|\begin{array}{c} \Xc_{15}\\ ----- \\
\Xc^*_{4,5}\end{array}\right|,$ \qquad\qquad $\pi(\Jc_{55})= \left|\begin{array}{ccccc}
s_{51}&s_{52}&s_{53}&s_{54}&s_{55 }\\
	s^*_{21}&s^*_{22}&s^*_{23}&s^*_{24}&s^*_{25 }\\
	s^*_{31}&s^*_{32}&s^*_{33}&0&0\\
	s^*_{41}&s^*_{42}&s^*_{43}&0&0\\
	s^*_{51}&0&0&0&0\\
\end{array}\right|$,\\
\\
$\Jc_{45} = \left|\begin{array}{c} \Xc_{25}\\ ----- \\ \Xc^*_{3,5}\end{array}\right|,$ 
\qquad\qquad $\pi(\Jc_{45})= \left|\begin{array}{ccccc}
s_{41}&s_{42}&s_{43}&s_{44}&s_{45 }\\
s_{51}&s_{52}&s_{53}&s_{54}&s_{55 }\\
s^*_{31}&s^*_{32}&s^*_{33}&0&0\\
s^*_{41}&s^*_{42}&s^*_{43}&0&0\\
s^*_{51}&0&0&0&0\\
\end{array}\right|$,\\
\\
$\Jc_{35} = \left|\begin{array}{c} \Xc_{35}\\ ----- \\ \Xc^*_{2,5}\end{array}\right|,$ 
\qquad\qquad $\pi(\Jc_{35})= \left|\begin{array}{ccccc}
0&0&s_{33}&s_{34}&s_{35 }\\
s_{41}&s_{42}&s_{43}&s_{44}&s_{45 }\\
s_{51}&s_{52}&s_{53}&s_{54}&s_{55 }\\
s^*_{41}&s^*_{42}&s^*_{43}&0&0\\
s^*_{51}&0&0&0&0\\
\end{array}\right|$,\\
\\
\\
$\Jc_{25} = \left|\begin{array}{c} \Xc_{45}\\ ----- \\ \Xc^*_{1,5}\end{array}\right|,$ 
\qquad\qquad $\pi(\Jc_{25})= \left|\begin{array}{ccccc}
0&0&s_{23}&s_{24}&s_{25 }\\
0&0&s_{33}&s_{34}&s_{35 }\\
s_{41}&s_{42}&s_{43}&s_{44}&s_{45 }\\
s_{51}&s_{52}&s_{53}&s_{54}&s_{55 }\\
s^*_{51}&0&0&0&0\\
\end{array}\right|$,\\
\\
$\Jc_{15} = \det \Xb$, ~~~ $\pi(\Jc_{15}) =\det S$.
}

\section{Fields of invariants for   radicals of  parabolic subgroups in the orthogonal and symplectic subgroups}

Let   $K$ be a field of   characteristic zero,  $G$ be the orthogonal or symplectic group over $K$,  and $P$ be its parabolic subgroup. There is a decomposition of the subgroup  $P$ into a product  $P=LU$ of a Levi subgroup  $L$  and the unipotent radical $U$. Let  $\Ac=K[G]$ be  the ring of regular functions on   $G$, and  $\Fc=K(G)$ be the field of rational functions on  $G$.
In this section, we construct a system of free generators  in the field of invariants  $\Fc^U$ with respect to the action of $U$ on $G$ by conjugation.  

We introduce necessary notations. For each matrix  $A$, we denote by $A^t$ (respectively, $A^\sigma$ ) its  transpose  (respectively, its transpose  with respect to the anti-diagonal). 
Let  $\Ib_N$ be the matrix of size $N\times N$ with ones on the anti-diagonal and zeros on other places. By definition, the orthogonal group  $\Om(N)$ is the group of all  $g\in\GL(N)$ such that  $g^t\Ib_N g=\Ib_N$.  

The symplectic group  $\Spm(N)$, where $N=2n$, is defined analogously by the matrix  $$\Jb_{2n}=\left(\begin{array}{cc}
0&-\Ib_n\\
\Ib_n&0
\end{array}\right).$$ 

Consider a decomposition of the integer segment $[1,N]$ into   a system of  consecutive segments $[1,n] =I_1\sqcup I_2\sqcup\ldots\sqcup I_\ell$. Denote   $n_k=|I_k|$. Suppose that this decomposition is symmetric with respect to the center of the segment, i.e. $I'_k=I_{\ell-k+1}$. 

Introduce notations  $\ell_0=\left[\frac{\ell}{2}\right]$, ~  $N_0=n_1+\ldots+n_{\ell_0}$. The cases of even and odd  $\ell$ separate.  
If $\ell=2\ell_0$, then  $N=2n$ and $N_0=n$. In the odd  case $\ell=2\ell_0+1$, denote $I_0=I_{\ell_0+1}$,~ $n_0=n_{\ell_0+1}$. Let $G_0$ stand for the group  $\Om(n_0)$ in the orthogonal case or 
the group $\Spm(n_0)$ in the symplectic case.

By this decomposition, we  construct the parabolic subgroup $P=LU$ with the Levi subgroup  $L$ of block-symmetric matrices  $\diag(A_1,\ldots,A_\ell)$,~
$A_k\in\GL(n_k)$, ~$A_{\ell-k+1}= (A_k^\sigma)^{-1}$ for each $1\leq k\leq \ell_0$ and  $A_{\ell_0+1}\in G_0$ (for $\ell=2\ell_0+1$).

If $\ell=2\ell_0$, then  each  $g\in P$ has the form  

$$g=\left(\begin{array}{cc}
A&0\\
0&(A^{\sigma})^{-1}
\end{array}\right)\left(\begin{array}{cc}
E&B\\
0&E
\end{array}\right) =
\left(\begin{array}{cc}
A&AB\\
0&(A^{\sigma})^{-1}
\end{array}\right),$$
where $A$ is a matrix from the parabolic subgroup in $\GL(N_0)$, which relates to the decomposition  $I_1\sqcup
 \ldots \sqcup I_{\ell_0}$ of the segment  $[1,N_0]$, and  $B$ is a skew-symmetric  (respectively, symmetric) matrix with respect to the anti-diagonal  in the orthogonal case
   (respectively, in the symplectic case).

 Consider the case  $\ell=2\ell_0+1$.  
   Denote  $\Ib_0=\Ib_{N_0}$,~ $$\Jb_0=\left\{\begin{array}{l}
 \Ib_{n_0},~~\mbox{in ~the ~case}~~ \Om(N),\\
  \Jb_{n_0}, ~~\mbox{in ~the ~case}~~ \Spm(N).
 \end{array}\right.$$

Each element  $g\in P$ has the form   
 $$g=\left(\begin{array}{ccc}
 A&0&0\\
 0&A_0&0\\
 0&0&(A^{\sigma})^{-1}
 \end{array}\right)
 \left(\begin{array}{ccc}
 E&V&\frac{1}{2}V W\\
 0&E&W\\
 0&0&E
 \end{array}\right) 
 \left(\begin{array}{ccc}
 E&0&B\\
 0&E&0\\
 0&0&E
 \end{array}\right) = $$
 $$
 \left(\begin{array}{ccc}
 A&AV& A\left(B+\frac{1}{2}V W \right)\\
 0&A_0&A_0 W\\
 0&0&(A^{\sigma})^{-1}
 \end{array}\right),$$  
where $A$ and $B$  as above, $V$ is a matrix of size  $N_0\times n_0$,~ $W=-\Jb_0V^t \Ib_0$,  and
$A_0\in G_0$.

 Denote $S^\circ=S\cap G$, where $S$ is the subset of $\GL(N)$ defined in the previous section.  Arguing analogously   Proposition \ref{open}, we obtain that 
  the subset   $$\bigcup_{g\in U} gS^\circ g^{-1}$$
 is dense in $G$. 
 
 The restriction map  $$\pi:\Ac \to K[S^\circ]$$ establishes the 
 embedding of the field of $U$-invariants $\Fc^U$ into  $K(S^\circ)$. 

If $\ell=2\ell_0$, then   $S^\circ$  consists of all matrices of the form 
$$\left(\begin{array}{cc}
0&S_{12}\\
S_{21}&S_{22}\\
\end{array}\right) = \left(\begin{array}{cc}
0& \pm\Ib_0 (A^{\sigma})^{-1}\\
\Ib_0A& \Ib_0 AB\\
\end{array}\right).$$

If $\ell=2\ell_0+1$,  then  $S^\circ$  consists of all matrices of the form 
\begin{equation}\label{Smatrix}
\left(\begin{array}{ccc}
0&0&S_{13}\\
0&S_{22}&S_{23}\\
S_{31}&S_{32}&S_{33}
\end{array}\right)= \left(\begin{array}{ccc}
0&0&  \pm\Ib_0 (A^{\sigma})^{-1}\\0&\Jb_0A_0& \Jb_0A_0W\\
\Ib_0 A &\Ib_0AV& \Ib_0 A\left(B+\frac{1}{2}V W\right)
\end{array}\right).
\end{equation}

Consider the  set  $\Sb^\circ $ of pairs  $(i,j)$ such that  \\
1) ~$i>N_0$ (for $\ell=2\ell_0$) or $i>N_0+n_0$ (for $\ell=2\ell_0+1$),\\
	2)~ if $i\in I_k$ and $j\in I'_m$, then $k\geq m$,\\
	3) ~ $i>j$  in the orthogonal case or $i\geq j$ in the symplectic case.

Similarly to the previous section, let  $\Sb^\circ_0$ consist of all pairs  $(i,j)\in \Sb^\circ$ lying on or above the anti-diagonal, respectively,  $\Sb^\circ_1$ consist of all pairs $(i,j)$ lying below the anti-diagonal.
Denote by $\Gb_0$  the set  $I_0\times I_0$.

\Ex~\Num\label{extwo}. Let $G=\Spm(8)$  and the parabolic subgroup is defined by the decomposition 
 $[1,8]=\{1\}\sqcup [2,3]\sqcup [4,5]\sqcup [6,7]\sqcup \{8\}$. On the   $(6\times 6)$-tables below, we mark by the symbol "$\times$"\,  the  cells $(i,j)$ belonging to $\Sb^\circ$ (respectively, belonging to  $\Sb^\circ_0$, ~$\Sb^\circ_1$,~ $\Gb_0$). We have
 {\scriptsize   $$\Sb^\circ=\left(\begin{array}{cccccccc}
\cdot&\cdot&\cdot&	\cdot&	\cdot&	\cdot&	\cdot&	\cdot\\
	\cdot&\cdot&\cdot&	\cdot&	\cdot&	\cdot&	\cdot&	\cdot\\
	\cdot&\cdot&\cdot&	\cdot&	\cdot&	\cdot&	\cdot&	\cdot\\
	\cdot&\cdot&\cdot&	\cdot&	\cdot&	\cdot&	\cdot&	\cdot\\
	\cdot&\cdot&\cdot&	\cdot&	\cdot&	\cdot&	\cdot&	\cdot\\
	\cdot&	\times&	\times&	\times&	 \times&\times&\cdot&	\cdot\\
			\cdot&	\times&	\times&		\times& \times&\times&\times&	\cdot\\
		\times&	\times&	\times&	\times&		\times&\times&\times&\times\\
			\end{array}\right),\qquad 		
				\Gb_0=\left(\begin{array}{cccccccc}
				\cdot&\cdot&\cdot&	\cdot&	\cdot&	\cdot&	\cdot&	\cdot\\
				\cdot&\cdot&\cdot&	\cdot&	\cdot&	\cdot&	\cdot&	\cdot\\
				\cdot&\cdot&\cdot&	\cdot&	\cdot&	\cdot&	\cdot&	\cdot\\
				\cdot&\cdot&\cdot&	\times&	\times&	\cdot&	\cdot&	\cdot\\
				\cdot&\cdot&\cdot&	\times&	\times&	\cdot&	\cdot&	\cdot\\
					\cdot&\cdot&\cdot&	\cdot&	\cdot&	\cdot&	\cdot&	\cdot\\	\cdot&\cdot&\cdot&	\cdot&	\cdot&	\cdot&	\cdot&	\cdot\\	\cdot&\cdot&\cdot&	\cdot&	\cdot&	\cdot&	\cdot&	\cdot\\
				\end{array}\right),$$}
			 {\scriptsize $$	\Sb^\circ_0=\left(\begin{array}{cccccccc}
				\cdot&\cdot&\cdot&	\cdot&	\cdot&	\cdot&	\cdot&	\cdot\\
				\cdot&\cdot&\cdot&	\cdot&	\cdot&	\cdot&	\cdot&	\cdot\\
				\cdot&\cdot&\cdot&	\cdot&	\cdot&	\cdot&	\cdot&	\cdot\\
				\cdot&\cdot&\cdot&	\cdot&	\cdot&	\cdot&	\cdot&	\cdot\\
				\cdot&\cdot&\cdot&	\cdot&	\cdot&	\cdot&	\cdot&	\cdot\\
				\cdot&	\times&\times&\cdot&\cdot&\cdot&\cdot&\cdot\\
					\cdot&	\times&\cdot&\cdot&\cdot&\cdot&\cdot&\cdot\\
						\times&	\cdot&\cdot&\cdot&\cdot&\cdot&\cdot&\cdot\\
				\end{array}\right),\qquad
		 \Sb^\circ_1=\left(\begin{array}{cccccccc}
				\cdot&\cdot&\cdot&	\cdot&	\cdot&	\cdot&	\cdot&	\cdot\\
				\cdot&\cdot&\cdot&	\cdot&	\cdot&	\cdot&	\cdot&	\cdot\\
				\cdot&\cdot&\cdot&	\cdot&	\cdot&	\cdot&	\cdot&	\cdot\\
				\cdot&\cdot&\cdot&	\cdot&	\cdot&	\cdot&	\cdot&	\cdot\\
				\cdot&\cdot&\cdot&	\cdot&	\cdot&	\cdot&	\cdot&	\cdot\\
				\cdot&	\cdot&	\cdot&	\times&	 \times&\times&\cdot&	\cdot\\
				\cdot&	\cdot&	\times&		\times& \times&\times&\times&	\cdot\\
				\cdot
				&	\times&	\times&	\times&		\times&\times&\times&\times\\
				\end{array}\right). 					$$}

	We attach to each  $(i,j)\in \Gb_0$  the rational function  $P_{ij}$ on  $G$  by the formula $$P_{ij}=\frac{M_{ij}}{M_0},$$ where
	$M_0$ is  the minor of order  $N_0$ with the systems of rows  $I_{\ell_0+2}\sqcup\ldots\sqcup I_\ell$ and columns $I_1\sqcup\ldots\sqcup I_{\ell_0}$, and  $M_{ij}$ is  the minor obtained from  $M_0$ by adding $i$th row and $j$th column.  
Easy to see that the minors  $M_{ij}$ and $M_0$ are invariants with respect to the action of  $U$ on  $G$ by conjugation. Therefore, $P_{ij}$ is also an invariant.

 Consider the subfield  $\Fc_0$ in $\Fc^U$ generated by the system  $\{P_{ij}:~(i,j)\in \Gb_0\}$.
	One can treat $G_0$ as a factor of the parabolic subgroup $P$.  This makes it possible to consider the field of rational functions  $K(G_0)$  as a subfield in the field  $K(G)$
  and also as a subfield in $K(S)$. 
	The restriction map   $\pi$ establishes an isomorphism of the subfield  $\Fc_0$ in  $\Fc^U$ onto the subfield   $K(G_0)$  in  $K(S)$.

For each pair  $(i,j)\in\Sb^\circ$, we denote by  $s_{ij}$ the restriction of the matrix element  $x_{ij}$ on $S^\circ$.\\
\\
\Lemma\Num\label{sss}. 
The field   $K(S^\circ)$ is freely generated over the subfield  $K(G_0)$  by the system of matrix elements 
$\{s_{ij}:~~ (i,j)\in\Sb^\circ\}$.\\
\Proof. Let us prove for   $\ell=2\ell_0+1$. The case  $\ell=2\ell_0$ is treated similarly.  
We consider the system of matrix elements  $$\{a_{ij}, v_{ij}, b_{ij}\}$$  
from (\ref{Smatrix}), where  $A=(a_{ij})$,~ $V=(v_{ij})$ and  $ \{b_{ij}\}$ are  entries of  $B$  lying above the anti-diagonal in the case of $\Om(N)$  (respectively, on or above the anti-diagonal in the case of  $\Spm(N)$).  We treat these matrix elements as rational functions on $S^\circ$.
 The formula  (\ref{Smatrix}) implies that the field  $K(S^\circ)$ is freely generated over the subfield  $K(G_0)$ by the system  $\{ a_{ij}, ~ v_{ij},~ b_{ij}\}$.
The subfield generated by the matrix elements from  $S_{31}$ and $S_{32}$ coincides with the subfield  generated by  $\{a_{ij},~ v_{ij}\}$.  Recall that   
$S_{33} = \Ib_0 A\left(B+\frac{1}{2}V W\right)$, where $A$ is a matrix from the parabolic subgroup  in  $\GL(N_0)$ defined by the decomposition  $I_1\sqcup
\ldots \sqcup I_{\ell_0}$ of the segment  $[1,N_0]$, ~~  $B$ is the skew-symmetric  (respectively, symmetric) matrix with respect to the anti-diagonal in the orthogonal case
(respectively, in the symplectic case), and  $W=-\Jb_0V^t \Ib_0$. Rewrite  
$S_{33}$ in the form
$$S_{33}= \Ib_0 A \Ib_0 \cdot C+ \frac{1}{2}\Ib_0 AV W,$$
where $C= \Ib_0B$  is a skew-symmetric matrix in the orthogonal case and symmetric in the symplectic case. 
Consider the system of matrix elements  $\{c_{ij}\}$ from  $C$ lying below the diagonal in the case  $\Om(N)$  (respectively, on or below the diagonal in the case of  $\Spm(N)$). 
The system  $\{c_{ij}\}$ coincides with the system of $\{b_{ij}\}$ defined above.  
Observe that  $\Ib_0 A \Ib_0$  is a  lower block-triangular matrix defined by the partition of the segment  $[1,N_0]$ into consecutive
segments of sizes  $n_{\ell_0},\ldots,n_1$. 
Present the matrix   $\Ib_0 A \Ib_0$ as a product of two matrices with rational entries  $\Ib_0 A \Ib_0= RQ$, where $Q$ is a lower triangular matrix and $R$ is an upper unitriangular matrix.
Then the elements   
$s_{ij}$,~ $(i,j)\in\Sb^\circ$ from  $S_{33}$ are obtained  from  $\{c_{ij}\}$ by a composition of two triangular transformations:
\begin{itemize}
	\item $C\to QC$  is a triangular transformation with respect to the order relation: 
	$(i_1,j_1)$ is less than  $(i,j)$ if $j_1<j$ or $j_1=j$,~ $i_1<i$,
\item ~$QC\to  R Q C+ \frac{1}{2}\Ib_0 AV W$ is a triangular transformation with respect to the order relation  $\prec$ from  (\ref{Order}).
\end{itemize}
The system  $\{s_{ij}:~~ (i,j)\in\Sb^\circ\}$  is obtained from  $\{ a_{ij}, ~ v_{ij},~ b_{ij}\}$ by the composition of two triangular transformations.~ $\Box$

Denote by  $\Jc_{ij}^{\circ}$ the restriction of the determinant $\Jc_{ij}$ on the orthogonal or symplectic group  $G$.\\
\Prop\Num.~~ $\Jc_{ij}^\circ\ne 0$ for each $(i,j)\in \Sb^\circ$.\\
\Proof. For  $(i,j)\in \Sb_0^\circ$ the statement is obvious.
Consider the case  $(i,j)\in \Sb_1^\circ$.   
 Define the lower triangular with respect to the anti-diagonal  matrix   $A=(a_{st})\in S^\circ$ that has  nonzero entries on the anti-diagonal, has non-degenerate   
   $(i',i')$-block defined by the system of last  $i'$ rows and  the columns  $[j-i'+1, j]$, and  $a_{st}=0$ on other places. Then the determinant  $\Jc_{ij}^\circ(A)$ has the form 
$$\left|\begin{array}{cc}
A_{11}&A_{12}\\
A_{21}&0
\end{array}\right|,$$
where $A_{12}$ is the defined above non-degenerate  $(i',i')$-block,  and $A_{21}$ is an upper triangular block with non-zeros on the anti-diagonal. 
 Therefore $\Jc_{ij}^\circ(A)\ne 0$.~ $\Box$ \\
\Cor\Num. ~ $\pi(\Jc^\circ_{ij})\ne 0$ for each  $(i,j)\in \Sb^\circ$. \\

Recall that the subfield $\Fc_0$ in  $\Fc^U$ is  generated by   $\{P_{ij}:~(i,j)\in \Gb_0\}$, and it is isomorphic to $K(G_0)$.\\
\Theorem\Num\label{theoremtwo}. For the orthogonal or symplectic group $G$, the  system of polynomials  $\{\Jc^\circ_{ij}: ~(i,j)\in \Sb^\circ\}$ freely generates over   $\Fc_0$ the field of invariants  $\Fc^{U}$. \\
\Proof. 
Consider the extension   $\Pc_0$ generated over $\Fc_0$  by   $$\{\pi(\Jc_{ij}): ~~(i,j)\in \Sb^\circ_0 \}.$$
Let $Q$ be the parabolic subgroup in  $G$ defined by the partition of the segment  $[1, N]$ into the system of  consecutive segments of sizes  $(1,\ldots, 1, n_0, 1,\ldots,1)$,  and $U(Q)$ be its  unipotent radical. The group  $U(Q)$ acts on  $S^\circ$ by right multi\-pli\-ca\-tion. 
Similar to the proof of Proposition \ref{fzero}, one can show that  $\pi$ established isomorphism of the subfield  $\Pc_0$ onto the subfield of invariants in   $K(S^\circ)$ under the right   multiplication by  $U(Q)$.  
It follows from the proof of Theorem  \ref{theoremone} that the system of polynomials  $\{\pi(\Jc_{ij}^\circ):~ (i,j)\in\Sb^\circ\}$ is obtained by  a triangular transformation from the system  $\{s_{ij}:~(i,j)\in\Sb^\circ\}$. Applying Lemma \ref{sss} we conclude the proof. ~$\Box$

\end{document}